\documentclass[11pt,lecno]{amsart}
\usepackage{indentfirst,amssymb,amsmath,amsthm}
\newtheorem{theo}{Theorem}[section]

\newtheorem{remark}{Remark}[section]
\newtheorem{corollary}{Corollary}[section]

\newcommand{\be}{\begin{equation}}
\newcommand{\ee}{\end{equation}}
\newcommand{\bea}{\begin{eqnarray}}
\newcommand{\eea}{\end{eqnarray}}
\newcommand{\beas}{\begin{eqnarray*}}
\newcommand{\eeas}{\end{eqnarray*}}

\begin{document}

\setcounter{page}{1} \setlength{\unitlength}{1mm}\baselineskip
.58cm \pagenumbering{arabic} \numberwithin{equation}{section}

\title[ ]
{Perfect fluid spacetimes and $k$-almost yamabe solitons}

\author[ ]
{ Krishnendu De, Uday Chand De and Aydin Gezer $^{*}$ }

\address
 { Department of Mathematics,
 Kabi Sukanta Mahavidyalaya,
The University of Burdwan.
Bhadreswar, P.O.-Angus, Hooghly,
Pin 712221, West Bengal, India. ORCID iD: https://orcid.org/0000-0001-6520-4520}
\email{krishnendu.de@outlook.in }

\address
{ Department of Pure Mathematics, University of Calcutta, West Bengal, India. ORCID iD: https://orcid.org/0000-0002-8990-4609}
\email {uc$_{-}$de@yahoo.com}

\address
{$^{*}$  Department of Mathematics, Ataturk University, Erzurum-TURKEY}
	\email {aydingzr@gmail.com}

\footnotetext {$\bf{2020\ Math.\ Sub.\ Classification\:}.$ 83D05; 83C05; 53C50.
\\ {Key words: Perfect fluids; k-almost Yamabe solitons; Robertson-Walker spacetimes .\\
\thanks{$^{*}$ Corresponding author}
}}
\maketitle

\vspace{1cm}

\begin{abstract}
In this article, we presumed that a perfect fluid is the source of the gravitational field while analyzing the solutions to the Einstein field equations. With this new and creative approach, here we study $k$-almost yamabe solitons and gradient $k$-almost yamabe solitons. First, two examples are constructed to ensure the existence of gradient $k$-almost Yamabe solitons. Then we show that if a perfect fluid spacetime admits a $k$-almost yamabe soliton, then its potential vector field is Killing if and only if the divergence of the potential vector field vanishes. Besides, we prove that if a perfect fluid spacetime permit a $k$-almost yamabe soliton ($g,k,\rho,\lambda$), then the integral curves of the vector field $\rho$ are geodesics, the spacetime becomes stationary and the isotopic pressure and energy density remain invariant under the velocity vector field $\rho$. Also, we establish that if the potential vector field is pointwise collinear with the velocity vector field and $\rho(a)=0$ where a is a scalar, then either the perfect fluid spacetime represents phantom era, or the potential function $\Phi$ is invariant under the velocity vector field $\rho$. Finally, we prove that if a perfect fluid spacetime permits a gradient $k$-almost yamabe soliton ($g,k,D\Phi,\lambda$) and $R, \lambda, k$ are invariant under $\rho$, then the vorticity of the fluid vanishes.

\end{abstract}

\maketitle

\section{Introduction}

In this article, we will deal with the spacetimes that obey the Einstein field equations (in short, EFEs) when a perfect fluid (in short, PF) serves as the source of the gravitational field. If a Lorentzian manifold's non-vanishing Ricci tensor $S$ fulfills the conditions
\begin{equation}
\label{1}
S=\alpha g+ \beta C \otimes C,
\end{equation}
it is referred to as a PF spacetime, in which $\alpha$, $\beta$ (not simultaneously zero) are scalars and for any $X_{1}$, $g(X_{1}, \rho)=C(X_{1})$, $\rho$ stands for the velocity vector field.\par

Let $(M^{n}, g)$ be a Lorentzian manifold whose metric $g$ is of signature \newline $(-,+, +, \ldots, + )$, that is, $g$ is of index 1. In \cite{alias1}, Alias et al. proposed the notion of generalized Robertson-Walker (in short, GRW) spacetimes. A GRW spacetime is a Lorentzian manifold $M^{n}$ $(n \ge 3)$ that may be expressed as $M=-I \times  f^2  M^{*}$, in which the open interval $I$ contained in $\mathbb{R}$, $ M^{*(n-1)}$ denotes the Riemannian manifold and $f>0$ is a smooth function, named as scale factor or warping function. The above stated spacetime turns into Robertson-Walker (in short, RW) spacetime when the dimension of $M^{*}$ is three and the sectional curvature is constant.\par


Hamilton \cite{hamilton} initially proposed the fascinating concept of Yamabe flow while simultaneously introducing the Ricci flow to address Yamabe's conjecture. According to Yamabe's conjecture, a compact connected Riemannian manifold and a manifold with constant scalar curvature are conformally equivalent. This theory was described by Yamabe in 1960 \cite{yam}. A partial differential equation of the heat kind is described by
\begin{equation}\nonumber
\frac{\partial}{\partial t}g = -Rg,\hspace{.2cm} g(0) = g_{\circ},
\end{equation}
in which $R$ stands for the scalar curvature. The foregoing equation is called a Yamabe flow on a Riemannian manifold.\par
If a Riemannian metric g of a Riemannian manifold fulfills
\begin{equation}\label{2}
  \pounds_Z g=(R-\lambda)g,
\end{equation}
it is referred to as a Yamabe soliton (in short, YS), in which $\lambda$, $Z$ and $\pounds$ stand for real numbers, smooth vector fields, and the Lie derivative, respectively.\par

Similar to the Riemannian context, YSs have been characterized in semi-Riemannian manifolds in \cite{calvi}.\par
In a latest study, Chen and Deshmukh \cite{chen} introduced the novel idea of quasi-YS which is a generalization of YS on a Riemannian manifold as
\begin{equation}\label{3}
  (\pounds_Z g)(X_{1},Y_{1})=(R-\lambda)g(X_{1},Y_{1})+\;\Psi Z^*(X_{1})\;Z^*(Y_{1}),
\end{equation}
where  $Z^*$, $\Psi$ and $\lambda$ denote, in that order, the dual 1-form of $Z$, the smooth function, and a real integer. It is referred to as a proper quasi-YS if $\Psi \neq 0$.\par

Pirhadi and Razavi looked into a gradient almost quasi-YS in \cite{pr} several years ago, treating $\lambda$ as a smooth function. They received a few intriguing formulae. Chen \cite{ch} has recently investigated almost quasi-YSs in relation to almost Cosymplectic manifolds. Additionally, Blaga has investigated this kind of soliton on warped products \cite{blaga}.\par

Chen \cite{ch1} introduced a novel idea known as the $k$-almost Yamabe soliton (in short, $k$-AYS) in a current paper. Chen claims that if a nonzero function k, a smooth vector field $Z$ and a smooth function $\lambda$  exist such that
\begin{equation}\label{4}
  \frac{k}{2} \pounds_Z g- (R-\lambda)g=0
\end{equation}
holds, then a Riemannian metric is a $k$-AYS. The k-AYS is designated by the symbol $(g,Z,k,\lambda)$. The foregoing equation turns into a gradient $k$-AYS ($g, \Phi, k, \lambda$) for some smooth function $\Phi$ if $Z=D\Phi$. Then the previous equation takes into the following form
 \begin{equation}\label{5}
  Hess \Phi= \frac{1}{k}(R-\lambda)g,
\end{equation}
in which $Hess$ denotes the Hessian. If $Z=0$, the $k$-AYS is trivial; otherwise, nontrivial. Additionally, the preceding equation yields the $k$-YS when $\lambda=$ constant.\par
Recently, in PF spacetimes, several researchers studied numerous type of solitons like YSs \cite{de}, gradient YSs\cite{dez}, Ricci solitons (\cite{blaga2},\cite{dms}), gradient Ricci solitons(\cite{dez}, \cite{dms}), Ricci-Yamabe solitons\cite{sing}, gradient $\eta$-Einstein solitons(\cite{dms}), gradient m-quasi Einstein solitons\cite{dez}, gradient Schouten solitons(\cite{dms}),  respectively.\par

The above mentioned studies tells us that many recent investigations on gradient solitons have been published on Lorentzian manifolds. In \cite{de}, De et al. have studied YSs in PF-spacetimes and in this paper we investigate on $k$-AYSs in PF-spacetimes which is a natural generalization of YSs. Specifically, we prove the following theorems:
\begin{theo}
\label{thm1}
If a PF spacetime admits a $k$-AYS, then its potential vector field is Killing if and only if $div Z=0$.
\end{theo}
\begin{theo}\label{thm2}
If a PF spacetime permits a $k$-AYS ($g,k,\rho,\lambda$), then \par
(i)the integral curves of the vector field $\rho$ are geodesics,\par
(ii) the spacetime becomes stationary and \par
(iii) the energy density $\sigma$ and isotopic pressure p remain invariant under 
 $\rho$.
\end{theo}
\begin{theo}\label{thm3}
Let a PF spacetime permit a $k$-AYS ($g,k,Z,\lambda$). If 
Z is pointwise collinear with 
$\rho$ and for a scalar $a$, $\rho(a)=0$, then either the PF spacetime represents phantom era, or  potential function $\Phi$ is invariant under 
$\rho$.
\end{theo}

\begin{theo}\label{thm4}
Let a PF spacetime permit a gradient $k$-AYS ($g,k,D\Phi,\lambda$). If $R, \lambda, k$ are invariant under $\rho$, then the vorticity of the fluid vanishes.
\end{theo}
\begin{remark}
  Since every RW-spacetime is a PF-spacetime \cite{neil}, all the results of PF-spacetimes which have established in this paper are also true in RW-spacetime.
\end{remark}

\section{ perfect fluid spacetimes}
The equation (\ref{1}) yields
\begin{equation}
\label{3a}
Q X_{1}=\alpha X_{1}+\beta C(X_{1}) \rho,
\end{equation}
in which $Q$ is the Ricci operator described by $g(QX_{1}, Y_{1})=S(X_{1},Y_{1})$.\par

The covariant derivative of (\ref{3a}) provides
\begin{equation}
\label{5b}
(\nabla_{X_{1}}Q)(Y_{1})=X_{1} (\alpha)Y_{1}+X_{1}(\beta) C(Y_{1})\rho+\beta (\nabla_{X_{1}}C)(Y_{1})\rho+\beta C(Y_{1}) \nabla_{X_{1}}\rho.
\end{equation}

In the absence of the cosmological constant, Einstein's field equations have the following structure
\begin{equation}
\label{1.2}
S-\frac{R}{2}g=\kappa T,
\end{equation}
if $\kappa$ stands for the gravitational constant and $T$ for the energy momentum tensor.\par

In a PF spacetime $T$ is described by
\begin{equation}
\label{1.1}
T=(p+\sigma)C \otimes C+p g,
\end{equation}
in which $p$ is the isotropic pressure and $\sigma$ denotes the energy density.
The equations (\ref{1}), (\ref{1.2}) and (\ref{1.1}) together yield
\begin{equation}
\label{1.7}
\alpha=\frac{\kappa(p-\sigma)}{2-n},\quad\beta=\kappa(p+\sigma).
\end{equation}

Additionally, an equation of state (briefly, EOS) with the shape $p = p(\sigma)$ connects $p$ and $\sigma$, and the PF-spacetime is known as isentropic. Furthermore, if $p = \sigma$, the PF-spacetime is referred to as stiff matter. The PF- spacetime is referred to as the dark matter era if $p+\sigma=0 $, the dust matter fluid if $p = 0$, and the radiation era if $p =\frac{\sigma}{3}$ \cite{cha1}. The universe is represented as accelerating phase when $\frac{p}{\sigma}< {-\frac{1}{3}}$. It covers the quintessence phase if $-1< \frac{p}{\sigma}< 0$ and phantom era if $\frac{p}{\sigma}< -1$.\par

\section{Examples of gradient $k$-almost yamabe solitons}
In order to establish the validity of gradient k-AYSs, we set two examples of spacetime. While the second example has a non-constant potential function, the first one is of 
constant potential function.\par
If the function $\Phi$ is smooth, we obtain
\begin{equation}\label{ee1}
    \Phi,_i =\frac{\partial \Phi}{\partial w_{i}}\;\;
    \Phi ,_{ij}=\frac{\partial^2 \Phi}{\partial w_{i}\partial w_{j}}-\Gamma^{k}_{ij}\Phi ,_i.
\end{equation}
The equation (\ref {5}) can be described as follows (in local coordinate system):
\begin{equation}\label{ee2}
    \Phi ,_{ij}=\frac{1}{k}(R-\lambda)g_{ij}.
\end{equation}

\subsection{Example 1:}
We choose a Lorentzian metric $g$, expressed by
\begin{equation}\label{ee.1}
    ds^{2}=g_{ij}dw^{i}dw^{j}=(dw_1)^{2}+ (w_1)^2 (dw_2)^2+(w_2)^2(dw_3)^{2}-(dw_4)^{2},
\end{equation}
in a 4-dimensional Lorentzian manifold $\mathbb{R}^4$, in which $i,j=1,2,3,4$.\par
Using (\ref{ee.1}), the Lorentzian metric's non-vanishing components are stated by
\begin{equation}\label{ee.2}
    g_{11}=1,\;\; g_{22}=(w_1)^2,\;\;\ g_{33}=(w_2)^2,\;\; g_{44}=-1.
\end{equation}

Using the equation (\ref{ee.2}), the components of the non-vanishing Christoffel symbols are described by:
$$\Gamma^{1}_{22}=-w_{1},\;\;\Gamma ^{2}_{12}=\frac{1}{w_1},\;\;\ \Gamma ^{2}_{33}=-\frac{w_2}{(w_1)^2},\;\;\Gamma ^{3}_{23}=\frac{1}{w_2}.$$
Let us choose $\Phi (w_1, w_2, w_3, w_4)$ on $M$, an arbitrary smooth function. Hence, the followings are calculated
\begin{equation*}
    \Phi ,_{11}=\frac{\partial^2 \Phi}{(\partial w_{1})^2}\;\;,
\end{equation*}
\begin{equation*}
    \Phi ,_{12}=\frac{\partial^2 \Phi}{\partial w_{1}\partial w_{2}}-\frac{1}{w_1}\frac{\partial \Phi}{\partial w_{2}}=\Phi , _{21}\;\;,
\end{equation*}
\begin{equation*}
    \Phi ,_{13}=\frac{\partial^2 \Phi}{\partial w_{1}\partial w_{3}}=\Phi , _{31}\;\;,
\end{equation*}
\begin{equation*}
    \Phi ,_{14}=\frac{\partial^2 \Phi}{\partial w_{1}\partial w_{4}}=\Phi , _{41}\;\;,
\end{equation*}
\begin{equation*}
    \Phi ,_{22}=\frac{\partial^2 \Phi}{(\partial w_{2})^2}-w_1 \frac{\partial \Phi}{\partial w_1}\;\;,
\end{equation*}
\begin{equation*}
    \Phi ,_{23}=\frac{\partial^2 \Phi}{\partial w_{2}\partial w_{3}}-\frac{1}{w_2}\frac{\partial \Phi}{\partial w_{3}}=\Phi , _{32}\;\;,
\end{equation*}
\begin{equation*}
    \Phi ,_{24}=\frac{\partial^2 \Phi}{\partial w_{2}\partial w_{4}}=\Phi , _{42}\;\;,
\end{equation*}
\begin{equation*}
    \Phi ,_{34}=\frac{\partial^2 \Phi}{\partial w_{3}\partial w_{4}}=\Phi , _{43}\;\;,
\end{equation*}
\begin{equation*}
    \Phi ,_{33}=\frac{\partial^2 \Phi}{(\partial w_{3})^2}-\frac{w_2}{(w_1)^2} \frac{\partial \Phi}{\partial w_2}\;\;,
\end{equation*}
\begin{equation*}
    \Phi ,_{44}=\frac{\partial^2 \Phi}{(\partial w_{4})^2}\;\;.
\end{equation*}
Using the aforementioned equations, we obtain from the equation(\ref{ee2})
\begin{equation*}
    \frac{\partial^2 \Phi}{(\partial w_{1})^2}=\frac{1}{k}(R-\lambda)\;\;,
\end{equation*}
\begin{equation*}
    \frac{\partial^2 \Phi}{\partial w_{1}\partial w_{2}}-\frac{1}{w_1}\frac{\partial \Phi}{\partial w_{2}}=0\;\;,
\end{equation*}
\begin{equation*}
    \frac{\partial^2 \Phi}{\partial w_{1}\partial w_{3}}=0\;\;,
\end{equation*}
\begin{equation*}
    \frac{\partial^2 \Phi}{\partial w_{1}\partial w_{4}}=0\;\;,
\end{equation*}
\begin{equation*}
    \frac{\partial^2 \Phi}{(\partial w_{2})^2}-w_1 \frac{\partial \Phi}{\partial w_1}=\frac{1}{k}(w_1)^2(R-\lambda)\;\;,
\end{equation*}
\begin{equation*}
    \frac{\partial^2 \Phi}{\partial w_{2}\partial w_{3}}-\frac{1}{w_2}\frac{\partial \Phi}{\partial w_{3}}=0\;\;,
\end{equation*}
\begin{equation*}
    \frac{\partial^2 \Phi}{\partial w_{2}\partial w_{4}}=0\;\;,
\end{equation*}
\begin{equation*}
    \frac{\partial^2 \Phi}{\partial w_{3}\partial w_{4}}=0\;\;,
\end{equation*}
\begin{equation*}
    \Phi ,_{33}=\frac{\partial^2 \Phi}{(\partial w_{3})^2}-\frac{w_2}{(w_1)^2} \frac{\partial \Phi}{\partial w_2}=\frac{1}{k}(w_2)^2 (R-\lambda)\;\;,
\end{equation*}
\begin{equation*}
    \frac{\partial^2 \Phi}{(\partial w_{4})^2}=-\frac{1}{k}(R-\lambda)\;\;.
\end{equation*}
According to the aforementioned equations, $R=\lambda$ and $\Phi$ should be a constant function.\par
As a result, the metric is a gradient $k$-AYS with a constant potential function.

\subsection{Example 2:}
 we choose a Lorentzian metric $g$, described by
\begin{equation}\label{eee.1}
    ds^{2}=g_{ij}dw^{i}dw^{j}=e^{w_1 +1}(dw_1)^{2}+ e^{w_1}[ (dw_2)^2+(dw_3)^{2}-(dw_4)^{2}],
\end{equation}
in a Lorentzian manifold $\mathbb{R}^{4}$, in which $i,j=1,2,3,4$.\par
Using (\ref{eee.1}), we find 
\begin{equation}\label{eee.2}
    g_{11}=e^{w_1 +1},\;\; g_{22}=e^{w_1},\;\;\ g_{33}=e^{w_1},\;\; g_{44}=-e^{w_1}.
\end{equation}

Using the equation (\ref{eee.2}), we acquire the components of the non-vanishing Christoffel symbols as:
$$\Gamma^{1}_{11}=\frac{1}{2},\;\;\Gamma ^{1}_{22}=\Gamma ^{1}_{33}=-\frac{1}{2e},\;\;\
\Gamma ^{1}_{44}=\frac{1}{2e},\;\;\Gamma ^{2}_{12}=\Gamma ^{3}_{13}=\Gamma ^{4}_{14}=\frac{1}{2}.$$
Let us set an arbitrary smooth function $\Phi (w_1, w_2, w_3, w_4)$ on $M$, and calculate the followings
\begin{equation*}
    \Phi ,_{11}=\frac{\partial^2 \Phi}{(\partial w_{1})^2}-\frac{1}{2}\frac{\partial \Phi}{\partial w_{1}}\;\;,
\end{equation*}
\begin{equation*}
    \Phi ,_{12}=\frac{\partial^2 \Phi}{\partial w_{1}\partial w_{2}}-\frac{1}{2}\frac{\partial \Phi}{\partial w_{2}}=\Phi , _{21}\;\;,
\end{equation*}
\begin{equation*}
    \Phi ,_{13}=\frac{\partial^2 \Phi}{\partial w_{1}\partial w_{3}}-\frac{1}{2}\frac{\partial \Phi}{\partial w_{3}}=\Phi , _{31}\;\;,
\end{equation*}
\begin{equation*}
    \Phi ,_{14}=\frac{\partial^2 \Phi}{\partial w_{1}\partial w_{4}}-\frac{1}{2}\frac{\partial \Phi}{\partial w_{4}}=\Phi , _{41}\;\;,
\end{equation*}
\begin{equation*}
    \Phi ,_{22}=\frac{\partial^2 \Phi}{(\partial w_{2})^2}+\frac{1}{2e} \frac{\partial \Phi}{\partial w_1}\;\;,
\end{equation*}
\begin{equation*}
    \Phi ,_{23}=\frac{\partial^2 \Phi}{\partial w_{2}\partial w_{3}}=\Phi , _{32}\;\;,
\end{equation*}
\begin{equation*}
    \Phi ,_{24}=\frac{\partial^2 \Phi}{\partial w_{2}\partial w_{4}}=\Phi , _{42}\;\;,
\end{equation*}
\begin{equation*}
    \Phi ,_{34}=\frac{\partial^2 \Phi}{\partial w_{3}\partial w_{4}}=\Phi , _{43}\;\;,
\end{equation*}
\begin{equation*}
    \Phi ,_{33}=\frac{\partial^2 \Phi}{(\partial w_{3})^2}+\frac{1}{2e} \frac{\partial \Phi}{\partial w_1}\;\;,
\end{equation*}
\begin{equation*}
    \Phi ,_{44}=\frac{\partial^2 \Phi}{(\partial w_{4})^2}+\frac{1}{2e} \frac{\partial \Phi}{\partial w_1}\;\;.
\end{equation*}
Using the previous formulae, we derive from (\ref{ee2})
\begin{equation*}
    \frac{\partial^2 \Phi}{(\partial w_{1})^2}-\frac{1}{2}\frac{\partial \Phi}{\partial w_{1}}=\frac{1}{k}e^{w_1 +1}(R-\lambda)\;\;,
\end{equation*}
\begin{equation*}
    \frac{\partial^2 \Phi}{\partial w_{1}\partial w_{2}}-\frac{1}{2}\frac{\partial \Phi}{\partial w_{2}}=0\;\;,
\end{equation*}
\begin{equation*}
   \frac{\partial^2 \Phi}{\partial w_{1}\partial w_{3}}-\frac{1}{2}\frac{\partial \Phi}{\partial w_{3}}=0\;\;,
\end{equation*}
\begin{equation*}
    \frac{\partial^2 \Phi}{\partial w_{1}\partial w_{4}}-\frac{1}{2}\frac{\partial \Phi}{\partial w_{4}}=0\;\;,
\end{equation*}
\begin{equation*}
   \frac{\partial^2 \Phi}{(\partial w_{2})^2}+\frac{1}{2e} \frac{\partial \Phi}{\partial w_1}=\frac{1}{k}e^{w_1}(R-\lambda)\;\;,
\end{equation*}
\begin{equation*}
    \frac{\partial^2 \Phi}{\partial w_{2}\partial w_{3}}=0\;\;,
\end{equation*}
\begin{equation*}
   \frac{\partial^2 \Phi}{\partial w_{2}\partial w_{4}}=0\;\;,
\end{equation*}
\begin{equation*}
   \frac{\partial^2 \Phi}{\partial w_{3}\partial w_{4}}=0\;\;,
\end{equation*}
\begin{equation*}
    \frac{\partial^2 \Phi}{(\partial w_{3})^2}+\frac{1}{2e} \frac{\partial \Phi}{\partial w_1}=\frac{1}{k}e^{w_1}(R-\lambda)\;\;,
\end{equation*}
\begin{equation*}
    \frac{\partial^2 \Phi}{(\partial w_{4})^2}+\frac{1}{2e} \frac{\partial \Phi}{\partial w_1}=\frac{1}{k}e^{w_1}(R-\lambda)\;\;.
\end{equation*}
We take $\frac{R-\lambda}{k}=$ constant and solving the foregoing equations, we acquire $\Phi =2c e^{w_1 +1}$, $c\in$ $\mathbb{R}$ .\par
Thus, the metric is a gradient $k$-AYS with a potential function $\Phi =2c e^{w_1 +1}$.

\section{Proof of the main Theorems}

{\bf Proof of the Theorem \ref{thm1} :}\par
Let the Lorentzian metric of the PF-spacetimes permit a $k$-AYS. Hence, we acquire
\begin{equation}\label{p1}
  k\; \pounds_{Z}g=2(R-\lambda)g.
\end{equation}
Using explicit form of the Lie derivative, we reveal
\begin{equation}\label{p2}
  k[g(\nabla_{X_{1}}Z,Y_{1})+g(X_{1},\nabla_{Y_{1}}Z)]=2(R-\lambda)g(X_{1},Y_{1}).
\end{equation}
Taking contraction of the foregoing equation yields
\begin{equation}\label{p3}
  k \;div Z= 8(R-\lambda),
\end{equation}
where div stands for the divergence.\par
Hence, equation (\ref{p1}) provides
\begin{equation}\label{p4}
  k\; \pounds_{Z}g=\frac{div Z}{4}g,
\end{equation}
where equation (\ref{p3}) is used.\par
From the above we conclude that $Z$ is Killing if and only if $div Z=0$.\par
Thus the proof is completed.

{\bf Proof of the Theorem \ref{thm2} :}\par
Let the $k$-AYS's potential vector field $Z=\rho$. Hence, equation (\ref{p2}) rewritten as
\begin{equation}\label{p5}
  k\; \pounds_{\rho}g=k[g(\nabla_{X_{1}}\rho,Y_{1})+g(X_{1},\nabla_{Y_{1}}\rho)]=2(R-\lambda)g(X_{1},Y_{1}).
\end{equation}
Since $\rho$ is unit timelike in a PF spacetime, hence, we acquire $g(\rho, \rho)=-1$. The covariant derivative of $g(\rho, \rho)=-1$ gives $g(\nabla_{X_{1}}\rho, \rho)=0$, for all $X_{1}$ and $\nabla$ is the Levi-Civita connection.\par
Now putting $Y_{1}=\rho$ in equation (\ref{p5}) and using the previous result, we obtain
\begin{equation}\label{p6}
  k\; \pounds_{\rho}\rho=2(R-\lambda)\rho.
\end{equation}
Also replacing $X_{1}=Y_{1}=\rho$ in (\ref{p5}) provides
\begin{equation}\label{p7}
  R=\lambda.
\end{equation}
Hence, equation (\ref{p6}) entails $\pounds_{\rho}\rho=0$. Also from (\ref{p5}), we obtain $\pounds_{\rho}g=0$. We know that a spacetime is stationary if it has a time-like Killing vector field (see \cite{dug}, p.73). Therefore, the spacetime becomes stationary.\par
 Also, if $\rho$ is Killing, then we have $\pounds_{\rho}p=0$ and $\pounds_{\rho}\sigma=0$ (see \cite{dug}, p.89).\par
Therefore, the proof is finished.

{\bf Proof of the Theorem \ref{thm3} :}\par
Let $Z=a\rho$, where $a$ is a smooth function, that is, $Z$ is pointwise collinear with $\rho$. Therefore, we have
\begin{equation}\label{p8}
  \nabla_{X_{1}}Z=X_{1}(a)\rho+a\nabla_{X_{1}}\rho
\end{equation}
Using the previous equation in (\ref{p2}), we acquire
\begin{equation}\label{p9}
  k[X_{1}(a)C(Y_{1})+ag(\nabla_{X_{1}}\rho,Y_{1})+Y_{1}(a)C(X_{1})+ag(X_{1},\nabla_{Y_{1}}\rho)]=2(R-\lambda)g(X_{1},Y_{1}).
\end{equation}
Putting $Y_{1}=\rho$ in the preceding equation yields
\begin{equation}\label{p10}
  k[-X_{1}(a)+\rho(a)C(X_{1})+ag(X_{1},\nabla_{\rho}\rho)]=2(R-\lambda)C(X_{1}).
\end{equation}
Also, replacing $X_{1}$ with $\rho$ gives
\begin{equation}\label{p11}
  k \rho (a)=(R-\lambda).
\end{equation}
Also contracting (\ref{p9}), we get
\begin{equation}\label{p12}
  k\; a\; div \rho=3(R-\lambda)
\end{equation}
Let $\rho(a)=0$, then equation (\ref{p11}) implies $R=\lambda$. Using this in equation (\ref{p12}) yields $div \rho=0$ which implies that the velocity vector field is conservative. The nature of a conservative vector field is always irrotational, thus we conclude that the PF has zero vorticity.\par
This ends the proof.\par
{\bf Proof of the Theorem \ref{thm4} :}\par
Consider a gradient $k$-AYS ($g,\Phi,k,\lambda$) on a PF spacetime. Then, equation (\ref{5}) is given by
\begin{equation}\label{d1}
k\nabla _{X_{1}}D\Phi=(R-\;\lambda)Y_{1}.
\end{equation}
Covariant derivative  of (\ref{d1}) yields
\begin{eqnarray}
  \label{d2}
k \nabla_{X_{1}}\nabla_{Y_{1}} D\Phi&=&X_{1} (R-\;\lambda)Y_{1}+(R-\;\lambda)\nabla_{X_{1}}Y_{1}\nonumber\\&&
-\frac{1}{k}X_{1}(k)(R-\;\lambda)Y_{1}.
\end{eqnarray}

Interchanging $X_{1}$ and $Y_{1}$ in (\ref{d2}) provides
\begin{eqnarray}
  \label{d3}
k\nabla_{Y_{1}}\nabla_{X_{1}}D\Phi&=&Y_{1} (R-\;\lambda)X_{1}+(R-\;\lambda)\nabla_{Y_{1}}X_{1}\nonumber\\&&
-\frac{1}{k}Y_{1}(k)(R-\;\lambda)X_{1}
\end{eqnarray}

and
\begin{equation}\label{d4}
  k\nabla_{[X_{1},\;Y_{1}]}\;D\Phi=(R-\lambda)[X_{1},\;Y_{1}].
\end{equation}
Utilizing (\ref{d1})-(\ref{d4}) and  together with $R(X_{1},\;Y_{1})V_{1} =\nabla_{X_{1}} \nabla_{Y_{1}}V_{1}-\nabla_{Y_{1}} \;\nabla_{X_{1}}V_{1}-\nabla_{[X_{1},Y_{1}]}V_{1}$, we acquire
\begin{eqnarray}
  k^2 R(X_{1},\;Y_{1})D\Phi&=& k [\{X_{1}(R-\lambda)\}Y_{1}]-k [\{Y_{1}(R-\lambda)\}X_{1}]\nonumber\\&&
-X_{1}(k)(R-\lambda)Y_{1}+Y_{1}(k)(R-\lambda)X_{1}.\label{d5}
\end{eqnarray}
Contracting the above equation, we get
\begin{equation}\label{d6}
  k^2 S(Y_{1}, D\Phi)=-3k[Y_{1}(R)-Y_{1}(\lambda)]-3 Y_{1}(k)(R-\lambda).
\end{equation}
Also, from (\ref{1}), we have
\begin{equation}\label{d7}
 S(Y_{1}, D\Phi)=\alpha Y_{1}(\Phi)+\beta \rho(\Phi)C(Y_{1}).
\end{equation}
From the previous two equations, we infer
\begin{equation}\label{d8}
  -\frac{3}{k}[Y_{1}(R)-Y_{1}(\lambda)]-\frac{3}{k^2} Y_{1}(k)(R-\lambda)=\alpha Y_{1}(\Phi)+\beta \rho(\Phi)C(Y_{1}).
\end{equation}
Putting $Y_{1}=\rho$ in the previous equation yields
\begin{equation}\label{d9}
  -\frac{3}{k}[\rho(R)-\rho(\lambda)]-\frac{3}{k^2} \rho(k)(R-\lambda)=(\alpha -\beta) \rho(\Phi).
\end{equation}
If $R,\lambda,k$ are invariant under $\rho$, then
\begin{equation}\label{d10}
  (\alpha-\beta)\rho(\Phi)=0.
\end{equation}
Hence, from the above we say that either $\alpha=\beta$, or $\alpha \neq \beta$.\par
Case i: If $\alpha=\beta$, then using (\ref{1.7}) we obtain
\begin{equation}\label{d11}
    \frac{p}{\sigma}=-\frac{1}{3},
\end{equation}
which implies that the PF spacetime represents the phantom era.\par
Case ii: If $\alpha \neq \beta$, then $\rho (\Phi)=0$, that is, $\Phi$ is invariant under 
$\rho$.\par

Hence, the theorem is proved.
\section{Discussion}
A particular class of solutions on which the metric changes through diffeomorphisms and dilation performs a significant role in the investigation of flow singularities, since they appear as probable singularity models. Solitons is a common term used to describe them.\par
GR is applied mathematics' greatest achievement. GR has long been thought as both the most difficult and elegant physics theory ever created. Understanding GR, which disregards quantum effects, is essential for comprehending cosmology. In GR theory, the universe's matter content is determined by selecting the appropriate EMT, which is acknowledged to behave as a PF-spacetime in cosmological models. Here, EFE completes the crucial step in the construction of the cosmological model. In a few fields, including plasma physics, astronomy, atomic physical science, and nuclear physics, PF-spacetime models in general relativity theory are of great importance.\par

In this current investigations, we construct two examples to ensure the existence of gradient $k$-AYSs. Then we demonstrate that a PF spacetime's potential vector field is killing if and only if the potential vector field's divergence vanishes if a $k$-AYS is admitted to it. 
Additionally, we demonstrate that under specific circumstances, the PF spacetime reflects a phantom era. Finally, we demonstrate that the fluid's vorticity vanishes when a PF spacetime allows a gradient $k$-AYS.

\section{Declarations}
\subsection{Ethical Approval} Not applicable
\subsection{Competing interests}The authors declare that they have no conflict of interest.
\subsection{Authors' contributions} All authors contributed equally to this work.
\subsection{Funding }Not applicable
\subsection{Availability of data and materials }Not applicable

\end{document}